\documentclass[10pt,english,a4paper]{article}
\usepackage{babel}
\usepackage{amssymb}
\usepackage{amsmath}
\usepackage{enumerate}
\begin{document}
\def\R{\mathbb R}
\def\Z{\mathbb Z}
\def\Co{\mathcal{C}}
\newtheorem{th-def}{Theorem-Definition}[section]
\newtheorem{theo}{Theorem}[section]
\newtheorem{lemm}{Lemma}[section]
\newtheorem{prop}{Proposition}[section]
\newtheorem{defi}{D\'efinition}[section]
\newtheorem{cor}{Corollary}[section]
\def\dem{\noindent \textbf{Proof: }}
\def\rem{\noindent \textbf{Remark: }}
\def\rems{\noindent \textbf{Remarks: }}
\def\fin{  $\square$}
\def\ep{\varepsilon}
\def\uep{u_{\varepsilon}}
\def\omp{\Omega_{+}}
\def\om-{\Omega_{-}}

\numberwithin{equation}{section}
{\LARGE\begin{center} Morrey-Campanato estimates for Helmholtz equations with two unbounded media
\end{center}}
{\large\begin{center} Elise Fouassier
\end{center}}
%\maketitle
{\small\begin{center}\textbf{Abstract}
\end{center}

We prove uniform Morrey-Campanato estimates for Helmholtz equations in the case of two unbounded inhomogeneous media separated by an interface. They imply weighted $L^2$-estimates for the solution. We also prove a uniform $L^2$-estimate {\em without weight} for the trace of the solution on the interface.\\

}

\section{Introduction}

We consider the following Helmholtz equation
\begin{equation}\label{h}
i\ep u+\Delta u+n(x)u=f(x),\ x\in\R^d.
\end{equation}
We assume that the refraction index $n$ is nonnegative and discontinuous at the interface between two unbounded inhomogeneous media $\omp$, $\om-$ such that $ \omp \cup\overline{\om-}=\overline{\omp}\cup\om-=\R^d$. We write 
$$
 n(x)=\begin{cases}
       n_{+}(x) & \text{if $x\in\omp$}\\
       n_{-}(x) & \text{if $x\in\om-$}.
\end{cases}
$$
We study here the limiting absorption principle, i.e.~the limit when $\ep\to 0^+$ of equation (\ref{h}): our goal is to prove bounds on $u$ that are uniform in $\ep$. A related question is the statement of the Sommerfeld radiation condition for the following equation 
\begin{equation}
\Delta u+n(x)u=f(x),\ x\in\R^d.
\end{equation}
This is an open problem in the full generality of our assumptions, stated below.\\

In this paper, we prove uniform Morrey-Campanato type estimates for this equation, using a multiplier method borrowed from \cite{pv1}. These bounds encode in the optimal way the decay: $|u(x)|\sim 1/|x|^{\frac{d-1}{2}}$ at infinity of the solution $u$. They imply weighted $L^2$-estimates for the solution $u$. We also prove a uniform $L^2$-estimate {\em without weight} for the trace of the solution on the interface, which states that $u$ carries essentially no energy on this set.\\

In order to state precisely our assumptions and results, we need the following notations. 

\noindent First, we use the Morrey-Campanato norm, defined for $u\in L^2_{loc}$,
\begin{equation}
\|u\|_{\dot{B}^*}=\sup_{R>0}\frac{1}{R}\int_{B(R)}|u(x)|^2dx,
\end{equation}
where $B(R)$ denotes the ball of radius $R$. We also use the following dual norm
\begin{equation}
\|f\|_{\dot{B}}=\sum_{j\in\Z}\Big(2^{j+1}\int_{C(j)}|f(x)|^2dx\Big)^{1/2},
\end{equation}
where $C(j)=\{x\in\R^d/2^j\leq |x|\leq 2^{j+1}\}$.\\
The duality is given by the easy estimate
\begin{equation}\label{dual}
\Big|\int fudx\Big|\leq \|f\|_{\dot{B}}\|u\|_{\dot{B}^*}.
\end{equation}
Note that both norms $\dot{B}$ and $\dot{B}^*$ are homogeneous in space.\\

We also denote the radial and tangential derivatives by
$$\frac{\partial}{\partial r}=\frac{x}{|x|}\cdot\nabla,\qquad \frac{\partial}{\partial \tau}=\nabla-\frac{x}{|x|}\frac{\partial}{\partial r}.$$
\indent On the other hand, we assume that the interface between the two media $\Gamma=\partial\omp=\partial\om-$ is a smooth surface (Lipschitz is enough). Let $d\gamma$ be the euclidian surface measure on $\Gamma$ and $\nu(x)$ be the unit normal vector at $x\in\Gamma$ directed from $\om-$ to $\omp$. We denote, for $x \in \Gamma$, the jump
$$[n](x)=n_{+}(x)-n_{-}(x).$$
Throughout this paper we will write $\nabla n$ instead of $\nabla n_{+}\textbf{1}_{\omp}+\nabla n_{-}\textbf{1}_{\om-}$, the derivative of $n$ outside the interface. Similarly, $\partial_d n$ will denote the partial derivative of $n$ with respect to the $x_d$ variable outside the interface.\\

The key assumptions we make on the interface $\Gamma$ and the refraction index $n$ are the following. We comment on these assumptions later.\\
(H1) There is $\alpha>0$ such that the $d$-th component of $\nu$ satisfies
\begin{equation}
\nu_d(x)\geq\alpha\quad \text{for all}\quad x\in\Gamma .
\end{equation}
(H2) $[n](x)$ has the same sign for all $x\in\Gamma$; the following notation will be convenient: let $\sigma=-$ if $[n]$ is non-negative and $\sigma=+$ if $[n]$ is non-positive.\\
(H3) $n\in L^{\infty}$, $n\geq 0$.\\
(H4)\begin{equation}
2\sum_{j\in\Z}\sup_{C(j)}\frac{(x\cdot \nabla n(x))_{-}}{n(x)}:=\beta_1 < \infty.
\end{equation}
(H5)\begin{equation}
\frac{1}{\alpha}\sum_{j\in\Z}\sup_{C(j)}2^{j+1}\frac{(\partial_{d} n(x))_{\sigma}}{n(x)}:=\beta_2 < \infty.
\end{equation}
(H6) $\beta_1+\beta_2<1$.\\

We are now ready to state our main theorem.
\begin{theo}\label{th}
For dimensions $d\geq 3$, assume (H1)-(H6). Then, the solution to the Helmholtz equation (\ref{h}) satisfies the following estimates :
\begin{multline}\label{estim}
\|\nabla u\|_{\dot{B}^*}^2+\|n^{1/2}u\|_{\dot{B}^*}^2+\|(x\cdot\nabla n)_{+}^{1/2} u\|_{\dot{B}^*}^2+\int_{\R^d}\frac{|\nabla_{\tau}u|^2}{|x|}dx+\sup_{R>0}\frac{1}{R^2}\int_{S(R)}| u|^2 d\sigma_R \\
+\int_{\Gamma}\big|[n]\big| \; |u|^2 d\gamma +\int_{\R^d}|\partial_d n|\;|u|^2dx
\leq C \|f\|_{\dot{B}}^2,
\end{multline}
where $C$ is a constant depending only on $\alpha$, $\beta_1$ and $\beta_2$.
\end{theo}

We would like to stress several aspects of our analysis. First, the homogeneity of the estimates and assumptions makes this theorem compatible with the high frequencies. The scaling invariance plays a fundamental role in the high frequency limit of Helmholtz equations (see Benamou \textit{et al}~\cite{bckp}, Castella \textit{et al}~\cite{cpr} for the case of a regular index of refraction). \\
\indent We would like to point out several terms in the left-hand side. The last two say that, in principle, the energy is not trapped at the interface and it mainly radiates in the directions where $\partial_d n$ vanishes. On the other hand, the first two terms essentially assert that $u$ and $\nabla u$ belong to the optimal space $\dot{B}^*$, provided $f\in\dot{B}$. This is the same Morrey-Campanato estimates as in the regular case. \\
\indent Let us now comment our assumptions on $n$. First, they allow some growth at infinity : $n$ {\em does not} go to a constant at infinity. The hypotheses (H2), (H4) and (H5) can be understood as conditions on the trajectories of the geometrical optics. The condition (H4) implies the dispersion of these trajectories. The conditions (H2) and (H5) ensure that the energy goes from one side of the interface to the other. The condition (H5) involves both the interface and the index: it becomes a weaker assumption on the index when the interface is close to a hyperplane ($\alpha\sim 1$). This type of assumptions is natural in the study of the high frequency limit where the link with Liouville's equations can be understood through Wigner transform (see L. Miller~\cite{mil} for a refraction result in the case of a sharp interface for Schrodinger equation, E. Fouassier~\cite{ef} for high frequency limit of Helmholtz equations with interface; for an account on high frequency limit for wave equations, see P.-L. Lions, T.Paul~\cite{lp}, and \cite{bckp}, \cite{cpr} for Helmholtz equations without interface).\\
\indent Taking $[n]=0$ in Theorem \ref{th}, i.e. the case without interface, our results boil down to the uniform estimate proved by B. Perthame, L. Vega~\cite{pv1} when the refraction index satisfies assumptions (H3)-(H4) with $\beta_1<1$. Under these assumptions, they also proved in~\cite{pv3}, \cite{pv2} an energy estimate saying that the energy $|u|^2$ mainly radiates in the directions of the critical points of $n_{\infty}$ (where $n_{\infty}$ is given by $ n(x)\to n_{\infty}(\frac{x}{|x|}),\ |x|\to\infty$). In our case, the corresponding energy estimate corresponds to the last two terms.\\
\indent In the case of two unbounded media, similar results, but not scaling invariant, were obtained in previous papers. Eidus~\cite{ei} first proved weighted-$L^2$-estimates, and the $L^2$-estimate on the trace of the solution on the interface, for a piecewise constant index of refraction. To do so, he assumed that the interface satisfied an extra "cone-like" shape condition, $|x\cdot\nu|\leq C$ for $x\in \Gamma$. Under the same assumption on the interface, Bo Zhang~\cite{z} also proved inhomogeneous $B^*$-$B$ estimates when $n$ is a long-range perturbation of a piecewise constant function. S. DeBievre, D.W. Pravica~\cite{deb} proved weighted-$L^2$-estimates in a very general context using Mourre's commutator method. In particular they considered the case of an index that is smooth outside of a compact set in the $x_d$-direction and only bounded in this compact set.\\
\indent Our proof is based on a multiplier method. Following B. Perthame, L. Vega~\cite{pv1}, we use a combination of a Morawetz-type multiplier and of an elliptic multiplier. Following Eidus~\cite{ei}, we combine it with a multiplier specific to the case with an interface for which one direction plays a particular role (the $x_d$-direction here). The first two multipliers allow us to control both $\nabla u$ and $u$ locally in $L^2$ by $\|f\|_{\dot{B}}$ and we estimate the integral over the interface of $|u|^2$ using the third multiplier. \\

The paper is organized as follows. In Sections 2 and 3, we present the basic multipliers and the particular choice we make here to prove Theorem \ref{th}. Then, in Section 4, we give another estimate containing the trace of $\nabla u$ on the interface.
 
\section{Basic identities}

\begin{lemm}
The solution to the Helmholtz equation (\ref{h}) satisfies the following four identities, for smooth real valued test functions $\varphi$, $\psi$,
\begin{multline}\label{eg1}
-\int_{\R^d}\varphi(x)|\nabla u(x)|^2+\frac{1}{2} \int_{\R^d}\Delta\varphi(x)|u(x)|^2+\int_{\R^d}\varphi(x)n(x)|u(x)|^2 \\
=\mathcal{R}e\int_{\R^d} f(x)\varphi(x) \bar{u}(x),
\end{multline}
\begin{equation}\label{eg2}
\ep \int_{\R^d}\varphi(x)|u(x)|^2-\mathcal{I}m\int_{\R^d}\nabla\varphi(x)\cdot\nabla u(x)\bar{u}(x) 
=\mathcal{I}m\int_{\R^d} f(x)\varphi(x) \bar{u}(x),
\end{equation}
\begin{multline}\label{eg3}
\int_{\R^d}\nabla \bar{u}(x)\cdot D^2\psi(x)\cdot\nabla u(x)-\frac{1}{4} \int_{\R^d}\Delta ^2\psi(x)|u(x)|^2
+\frac{1}{2}\int_{\R^d}\nabla n(x)\cdot\nabla\psi(x)|u(x)|^2 \\
-\ep \mathcal{I}m\int_{\R^d}\nabla\psi(x)\cdot\nabla u(x)\bar{u}(x) +\frac{1}{2}\int_{\Gamma}[n]\nu(x)\cdot\nabla\psi(x)|u(x)|^2 d\gamma(x)\\
=-\mathcal{R}e\int_{\R^d} f(x)(\nabla\psi(x)\cdot\nabla \bar{u}(x)+\frac{1}{2}\Delta\psi(x)\bar{u}(x)),
\end{multline}
\begin{multline}\label{eg4}
\frac{1}{2}\int_{\Gamma}[n]\nu_d(x)|u(x)|^2 d\gamma(x)+\frac{1}{2} \int_{\R^d}\partial_{d} n(x)|u(x)|^2 \\
=-\mathcal{R}e\int_{\R^d} f(x)\partial_d\bar{u}(x)-\ep\mathcal{I}m\int_{\R^d}u(x)\partial_{d}\bar{u}(x).
\end{multline}
\end{lemm}

\dem The identities (\ref{eg1}) and (\ref{eg2}) are obtained by multiplying the Helmholtz equation (\ref{h}) by $\varphi \bar{u}$ and then taking the real and imaginary parts. The identity (\ref{eg3}) is obtained using the Morawetz-type multiplier $\nabla\psi(x)\cdot\nabla \bar{u}(x)+\frac{1}{2}\Delta\psi(x)\bar{u}(x)$ and taking the real part. To get the last identity (\ref{eg4}), we use the multiplier $\partial_{d} \bar{u}$ and take the real part.
\fin

\begin{lemm}
The solution to the Helmholtz equation (\ref{h}) satisfies the following estimate
\begin{multline}\label{trace}
\int_{\Gamma}|[n]||u(x)|^2 d\gamma(x)+\frac{1}{\alpha}\int_{\R^d}(\partial_{d} n)_{\sigma}|u|^2 \\
\leq 
2\int_{\R^d} |f\partial_{d}\bar{u}|
+2\ep\big|\mathcal{I}m\int_{\R^d}u\partial_{d}\bar{u}\big|
+\frac{1}{\alpha}\int_{\R^d}(\partial_{d} n)_{-\sigma}|u|^2.
\end{multline}
\end{lemm}

\dem
This lemma follows directly from the identity (\ref{eg4}) using the hypothesis (H1).\fin

\section{Proof of theorem \ref{th}}
\indent The following proof in which all the details are included for the convenience of the reader is essentially adapted from \cite{pv1}, apart from the treatment of the terms involving the interface and the partial derivative with respect to the $x_d$-direction, which is the key difficulty of this paper.\\
\indent We derive the proof of theorem \ref{th} from the above identities. We make the following choice of test functions $\psi$ and $\varphi$, for $R>0$,

$$  \nabla\psi(x)=\begin{cases}
       x/R & \text{for}\ |x|\leq R\\\
       x/|x| & \text{for}\ |x|>R,
\end{cases}
$$
$$  \varphi(x)=\begin{cases}
       1/2R & \text{for}\ |x|\leq R\\\
        0 & \text{for}\ |x|>R.
\end{cases}
$$

\noindent We also need the following calculations (in the distributional sense)
$$  D^2_{ij}\psi(x)=\begin{cases}
       \delta_{ij}/R & \text{for}\ |x|\leq R\\\
        (\delta_{ij}|x|^2-x_ix_j)/|x|^3 & \text{for}\ |x|>R,
\end{cases}
$$
$$  \Delta\psi(x)=\begin{cases}
       d/R & \text{for}\ |x|\leq R\\\
       (d-1)/|x| & \text{for}\ |x|>R,
\end{cases}
$$

\noindent and the inequality
\begin{equation}\label{ineg}
\frac{1}{4}\int_{\R^d} v\Delta (2\varphi-\Delta\psi)\geq\frac{d-1}{4R^2}\int_{S(R)}vd\sigma_R
\quad\text{for}\ v\geq 0.
\end{equation}

We add the identity (\ref{eg1}) to (\ref{eg3}), which gives, using the inequality (\ref{ineg}), for the previous choice of $\psi$, $\varphi$,

\begin{align*}
\frac{1}{2R}\int_{B(R)}|\nabla u(x)|^2 +\frac{1}{2R}&\int_{B(R)}n|u(x)|^2
+\int_{|x|>R}\frac{1}{|x|}\bigg(|\nabla u(x)|^2-\bigg|\frac{x}{|x|}\cdot\nabla u\bigg|^2\bigg) \\
&+\frac{d-1}{4R^2}\int_{S(R)}|u|^2 d\sigma_R 
+\frac{1}{2} \int_{\R ^d}(\nabla\psi(x)\cdot\nabla n(x))_{+} |u(x)|^2 \\
\leq 
C\int_{\R^d} |f(x)|\frac{|u(x)|}{|x|}+C&\int_{\R^d} |f(x)||\nabla u(x)|+\frac{1}{2} \int_{\R^d}(\nabla\psi(x)\cdot\nabla n(x))_{-} |u(x)|^2 \\
&+C\ep\int_{\R^d}|u(x)||\nabla u(x)|
+\frac{1}{2}\int_{\Gamma}\big|[n]\big|\;|\nu\cdot\nabla\psi ||u(x)|^2 d\gamma.
\end{align*}

Then, we use the inequality (\ref{trace}), together with the bound $|\nabla\psi|\leq 1$ to estimate the trace term

\begin{align}
\frac{1}{R} \int_{B(R)}|\nabla u(x)|^2 +\frac{1}{R}&\int_{B(R)}|u(x)|^2
+\frac{1}{R} \int_{B(R)}(x\cdot\nabla n(x))_{+} |u(x)|^2 \\
&+\int_{|x|>R}\frac{|\nabla_{\tau} u(x)|^2}{|x|} 
+\frac{d-1}{2R^2}\int_{S(R)}|u|^2 d\sigma_R \\
\leq 
C\int_{\R^d} |f(x)|\frac{|u(x)|}{|x|}&+C\int_{\R^d} |f(x)||\nabla u(x)|
+C\ep\int_{\R^d}|u(x)||\nabla u(x)|\\
&+\int_{\R^d}(\nabla\psi(x)\cdot\nabla n(x))_{-} |u(x)|^2 
+\frac{1}{\alpha}\int_{\R^d}(\partial_d n(x))_{\sigma}|u(x)|^2.
\end{align}

Our task is to estimate the terms in the right-hand side of the last inequality. We separate them in three types : those containing the source $f$, those containing $\ep$ and those containing the index $n$.\\
We begin by the two terms containing $f$. Using the duality estimate (\ref{dual}), we get for all $\delta>0$
\begin{eqnarray*}
\int_{\R^d} |f||\nabla u|&\leq &\|\nabla u\|_{\dot{B}^*} \|f\|_{\dot{B}}\\
                               &\leq & \delta\|\nabla u\|_{\dot{B}^*}^2+C_{\delta}\|f\|_{\dot{B}}^2.
\end{eqnarray*}

For the second term, we use a Cauchy-Schwarz inequality to obtain, for all $\delta>0$
\begin{eqnarray*}
\int_{\R^d} |f(x)|\frac{|u(x)|}{|x|}
&\leq & \sum_{j\in\Z}\bigg(2^{-j}\int_{C(j)}\frac{|u|^2}{|x|^2}\bigg)^{1/2}
\bigg(2^{j}\int_{C(j)}|f|^2\bigg)^{1/2} \\
&\leq & \bigg(\sup_{R>0}\frac{1}{R^2}\int_{S(R)}|u|^2d\sigma_R \bigg)^{1/2}
\sum_{j\in\Z}\bigg(2^{j}\int_{C(j)}|f|^2\bigg)^{1/2} \\
&\leq & \delta\sup_{R>0}\frac{1}{R^2}\int_{S(R)}|u|^2d\sigma_R+C_{\delta}N(f)^2.
\end{eqnarray*}

Next, we turn to the term containing $\ep$. First, using the identities (\ref{eg2}) and (\ref{eg1}) with $\varphi=1$, one can notice
\begin{eqnarray*}
\ep\int_{\R^d}|u|^2 &\leq & \int_{\R^d}|f\bar{u}|,\\
\ep\int_{\R^d}|\nabla u|^2 &\leq & \|n\|_{L^{\infty}}\int_{\R^d}|f\bar{u}|.
\end{eqnarray*}
Hence, using again the duality estimate (\ref{dual}), we obtain for all $\delta>0$
\begin{eqnarray*}
\ep\int_{\R^d}|u||\nabla u| &\leq & C\| u\|_{\dot{B}^*} \|f\|_{\dot{B}}\\
 &\leq & \delta\|\nabla u\|_{\dot{B}^*}^2+C_{\delta}\|f\|_{\dot{B}}^2.
\end{eqnarray*}

The terms with $n$ remain. We write
\begin{eqnarray*}
 \int_{\R^d}(\nabla\psi\cdot\nabla n)_{-} |u|^2 
&\leq & \sum_{j\in\Z}\int_{C(j)}n|u|^2\frac{(x\cdot\nabla n(x))_{-} }{n|x|}\\
&\leq & \bigg(\sup_{R>0}\frac{1}{R}\int_{B(R)}|u|^2 \bigg)
\sum_{j\in\Z}\sup_{C(j)}2^{j+1}\frac{(x\cdot\nabla n)_{-} }{n|x|} .
\end{eqnarray*}
Hence,
$$
\int_{\R^d}(\nabla\psi(x)\cdot\nabla n(x))_{-} |u(x)|^2 \leq  \beta_1 \|u\|_{\dot{B}^*}^2.
$$
Similarly, we get 
$$ 
\frac{1}{\alpha}\int_{\R^d}(\partial_d n(x))_{\sigma}|u(x)|^2
\leq \beta_2 \|u\|_{\dot{B}^*}^2.
$$

Putting all these estimates together gives for all $\delta>0$
\begin{multline*}
\frac{1}{R} \int_{B(R)}|\nabla u(x)|^2 +\frac{1}{R}\int_{B(R)}|u(x)|^2
+\frac{1}{R} \int_{B(R)}(x\cdot\nabla n(x))_{+} |u(x)|^2 \\
+\int_{|x|>R}\frac{|\nabla_{\tau} u(x)|^2}{|x|} 
+\frac{d-1}{2R^2}\int_{S(R)}|u|^2 d\sigma_R \\
\leq 
(\beta_1+\beta_2+\delta)\|u\|_{\dot{B}^*}^2+C_{\delta}\|f\|_{\dot{B}}^2.
\end{multline*}

Hence, choosing $\delta$ small enough (depending on $\beta_1+\beta_2$) and taking the supremum with respect to $R$, we obtain 
\begin{eqnarray}\label{intermediaire}
&&\|\nabla u\|_{\dot{B}^*}^2+\|n^{1/2}u\|_{\dot{B}^*}^2+\int_{\R^d}\frac{|\nabla_{\tau}u|^2}{|x|}+
\|(x\cdot\nabla n)_{+}^{1/2} u\|_{\dot{B}^*}^2 \\
&&\qquad\qquad\qquad\qquad +\sup_{R>0}\frac{1}{R^2}\int_{S(R)}| u|^2 d\sigma_R 
 \leq C\|f\|_{\dot{B}}^2.
\end{eqnarray}

To end up the proof, it remains to estimate the two last terms in the left-hand side of (\ref{estim}). We use the inequality (\ref{trace}) and the previous bounds 
\begin{eqnarray*}
&&\int_{\Gamma}\big|[n]\big|\;|u(x)|^2 d\gamma+\int_{\R^d}(\partial_d n)_{-\sigma}|u|^2 \\
&&\leq  
C\int_{\R^d} |f\partial_d\bar{u}|
+C\ep\int_{\R^d}|u|\nabla u|
+C\frac{1}{\alpha}\int_{\R^d}(\partial_d n)_{\sigma}|u|^2 \\
&&\leq  C (\|f\|_{\dot{B}}^2+\|u\|_{\dot{B}^*}^2). 
\end{eqnarray*}
Using (\ref{intermediaire}) we obtain
$$
\int_{\Gamma}\big|[n]\big|\;|u(x)|^2 d\gamma+\int_{\R^d}(\partial_d n)_{-\sigma}|u|^2 
\leq  C \|f\|_{\dot{B}}^2.
$$
This ends the proof.

\section{Another trace estimate}
In this section, we assume that the interface $\Gamma$ is a hyperplane $\Gamma=\{x_d=0\}$. Then we prove the following extra uniform estimate 

\begin{theo}\label{th2}
Under the assumptions of Theorem \ref{th}, if we assume moreover that $\Gamma=\{x_d=0\}$ and there exists $\beta>0$ such that $\langle x\rangle^{1+\beta}|\nabla_{x'} n|\in L^{\infty}$, then
\begin{equation}
\int_{\Gamma}\big|[n]\big|\;|\nabla u(x)|^2 dx' \leq C \Big(\|f\|_{\dot{B}}^2+\|\nabla_{x'}f\|_{\dot{B}}^2\Big)
\end{equation}
where $C$ is a constant depending only on $\alpha$, $\beta_1$, $\beta_2$ and $\|\langle x\rangle^{1+\beta}\nabla_{x'} n\|_{L^{\infty}}$.
\end{theo}

To prove this theorem, we will need the following identity
\begin{lemm}
The solution to the Helmholtz equation (\ref{h}) satisfies
\begin{multline}\label{deux}
\frac{1}{2}\int_{\Gamma}[n]\;|\nabla_{x'}u(x)|^2 dx'-\frac{1}{2}\int_{\Gamma}[n]\;|\partial_d u(x)|^2 dx'-\frac{1}{2}\int_{\Gamma}[n^2]\;|u(x)|^2 dx'\\
-\int_{\R^d}\partial_d(n^2(x))|u|^2-\ep \mathcal{I}m\int_{\R^d}n(x) u(x)\partial_d\bar{u}(x) \\
=\mathcal{R}e\int_{\R^d} n(x)f(x)\partial_d\bar{u}(x).
\end{multline}
\end{lemm}
\dem
This identity is obtained by multiplying the equation (\ref{h}) by $n\partial_d\bar{u}$ and taking the real part.\fin

\noindent \textbf{Proof of Theorem \ref{th2}}:\\
Firstly, note that $\nabla_{x'}n\ u\in \dot{B} $: 
\begin{eqnarray*}
\|\nabla_{x'}n\ u\|_{\dot{B}}&\leq &\sum_{j\in\Z}\big(2^{j+1}\int_{C(j)}
\langle x\rangle^{-2-2\beta} |u|^2\big)^{1/2} \\
   &\leq & \sum_{j<0}\big(2^{j+1}\int_{C(j)}
\langle x\rangle^{-2-2\beta} |u|^2 \big)^{1/2}\\
    &&+\sum_{j\geq 0}\big(2^{j+1}\int_{C(j)}
\langle x\rangle^{-2-2\beta} |u|^2\big)^{1/2} 
\end{eqnarray*}
but, using (\ref{estim}), we have, for $\delta>0$,
$$ \|\langle x\rangle^{-\frac{1}{2}-\delta} u \|_{L^2}\leq C \|u\|_{\dot{B}^*} \leq  C \|f\|_{\dot{B}},
$$
hence
\begin{eqnarray*}
\|\nabla_{x'}n\ u\|_{\dot{B}}
   &\leq & C\|f\|_{\dot{B}}\sum_{j<0}2^{\frac{j+1}{2}}+C\sum_{j\geq 0}2^{-j\beta/2}
\|\langle x\rangle^{-\frac{1}{2}-\frac{\beta}{2}} u\|_{L^2}\\
  &\leq &  C\|f\|_{\dot{B}}.
\end{eqnarray*}
So, we can apply the estimate (\ref{estim}) to $\nabla_{x'}u$. Indeed, $\nabla_{x'}u$ satisfies the following Helmholtz equation
\begin{equation}
i\ep \nabla_{x'}u+\Delta \nabla_{x'}u+n(x)\nabla_{x'}u=\nabla_{x'}f(x)+\nabla_{x'}n\ u,
\end{equation}
so $\nabla_{x'}u$ satisfies

$$\int_{\Gamma}\big|[n]\big|\;|\nabla_{x'} u|^2 dx' \leq C \Big(\|(\nabla_{x'}f\|_{\dot{B}}^2 +\|f\|_{\dot{B}}\Big).
$$

\noindent Then, we use (\ref{deux}) to get
\begin{multline}\label{to}
\int_{\Gamma}\big|[n]\big|\;|\partial_d u|^2 dx' \leq 
\int_{\Gamma}\big|[n]\big|\;|\nabla_{x'}u|^2 dx'+\int_{\Gamma}\big|[n^2]\big|\;|u|^2 dx'
+\int_{\R^d}|\partial_d(n^2)||u|^2 \\
+\ep\int_{\R^d}n |u\partial_d\bar{u}|
+\int_{\R^d} n|f\partial_d\bar{u}|.
\end{multline}

\noindent Since $n\in L^{\infty}$, we have 
$$ \int_{\Gamma}\big|[n^2]\big|\;|u|^2 dx' \leq 2\|n\|_{\infty}\int_{\Gamma}|[n]||u|^2 dx'\leq  C\|f\|_{\dot{B}}^2 $$
and
$$\int_{\R^d}|\partial_d(n^2)||u|^2 \leq 2\|n\|_{\infty}\int_{\R^d}|\partial_d n||u|^2
\leq C \|f\|_{\dot{B}}^2.$$
\noindent Hence, estimating the two last terms in the rigth-hand side of (\ref{to}) as before and using (\ref{estim}), we get
$$\int_{\Gamma}\big|[n]\big|\;|\partial_d u|^2 dx' \leq C \Big(\|f\|_{\dot{B}}^2+\|\nabla_{x'}f\|_{\dot{B}}\Big),$$
and thus the theorem is proved.\\

\noindent\textbf{Acknowledgement} I would like to thank Luis Vega for valuable comments on this work.\\
This work has been partially supported by the "ACI Jeunes Chercheurs - \\ M\'ethodes haute fr\'equence pour les \'equations diff\'erentielles ordinaires, et aux d\'eriv\'ees partielles. Applications", by the GDR "Amplitude Equations and Qualitative Properties" (GDR CNRS 2103 : EAPQ) and the European Program 'Improving the Human Potential' in the framework of the 'HYKE' network HPRN-CT-2002-00282.

\end{document}